\documentclass[11pt]{article}
\usepackage{times,mathptm}
\usepackage{amsfonts,amscd,amssymb}
\newtheorem{definition}{Definition}[section]
\newtheorem{lemma}[definition]{Lemma}
\newtheorem{proposition}[definition]{Proposition}

\newtheorem{remarks}[definition]{Remarks}

\def\yd{\mbox{$_H^H{\cal YD}$}}
\def\va{\varepsilon}
\def\v{\varphi}

\def\tr{\triangleright}
\def\rh{\rightharpoonup}
\def\lh{\leftharpoonup}

\def\ra{\rightarrow}
\def\a{\alpha}
\def\b{\beta}

\def\l{\lambda}

\def\cd{\cdot}

\def\O{\Omega}
\def\ov{\overline}
\def\un{\underline}
\def\bea{\begin{eqnarray*}}
\def\eea{\end{eqnarray*}}
\def\beq{\begin{equation}}
\def\eeq{\end{equation}}

\newcommand{\smi}{\mbox{$S^{-1}$}}

\def\rawo\lonra{\longrightarrow}

\def\ot{\otimes}

\def\gggg{{\mathfrak g}}
\def\complex{{\mathbb C}}
\def\nrat{{\mathbb Q}}

\newcommand{\selabel}[1]{\label{se:#1}}
\newcommand{\seref}[1]{Section~\ref{se:#1}}
\newcommand{\lelabel}[1]{\label{le:#1}}

\newcommand{\prlabel}[1]{\label{pr:#1}}
\newcommand{\prref}[1]{Proposition~\ref{pr:#1}}

\newcommand{\eqlabel}[1]{\label{eq:#1}}
\newcommand{\eqref}[1]{(\ref{eq:#1})}

\newenvironment{proof}{{\it Proof.}}{\hfill $ \square $ \vskip 4mm}
\begin{document}
\title{The quantum double for quasitriangular quasi-Hopf
algebras\thanks{Research supported by the bilateral project
``Hopf Algebras in Algebra, Topology, Geometry and Physics" of the Flemish and
Romanian governments.}}
\author{D. Bulacu\thanks{This paper was written while the first author was
visiting
the Free University of Brussels, VUB (Belgium); he would like to
thank VUB for its warm hospitality.}\\
Faculty of Mathematics\\ University of Bucharest\\
RO-70109 Bucharest 1, Romania\and
S. Caenepeel
\\ Faculty of Applied Sciences\\
Free University of Brussels, VUB\\ B-1050 Brussels, Belgium}
\date{}
\maketitle

\begin{abstract}
Let $D(H)$ be the quantum double associated to a finite
dimensional quasi-Hopf algebra $H$, as in \cite{hn1} and
\cite{hn2}. In this note, we first generalize a result of Majid \cite{m3}
for Hopf algebras, and then prove that the quantum double of
a finite dimensional
quasitriangular quasi-Hopf algebra is a biproduct in the sense of \cite{bn}.
\end{abstract}

\section{Introduction}\selabel{0}
Drinfeld \cite{d} introduced the quantum double $D(H)$ of a finite
dimensional Hopf algebra $H$. $D(H)$ is a quasitriangular Hopf
algebra, containing $H$ and $H^{*{\rm cop}}$ as subcoalgebras, and
Drinfeld used the construction to find the $R$-matrix $R_h$ of
$U_h(\gggg)$, the Drinfeld-Jimbo quantum enveloping algebra
associated to a complex semisimple Lie algebra $\gggg$ - the index
$h$ means that we work over the complex power series
$\complex[[X]]$, see \cite{k} for more detail. In fact he gives an
explicit formula for $R_h$ in the case where $\gggg={\rm sl}(2)$;
explicit descriptions in the general case are given in \cite{kr}
and \cite{ls}. Several alternative descriptions of the quantum
double have appeared in the literature, but the most striking one
is perhaps Majid's double cross product construction \cite{maj}.\\
Quasi-bialgebras and quasi-Hopf algebras were introduced by
Drinfeld in \cite{d1}, in connection with the
Knizhnik-Zamolodchikov system of partial differential equations,
cf. \cite{k}. Quasi-Hopf algebras also appear in algebraic number
theory, bearing information on the absolute Galois group ${\rm
Gal}(\ov{\nrat}/\nrat)$, cf. \cite{d3}. Quasi-bialgebras and
quasi-Hopf algebras naturally come out of categorical
considerations, putting some additional structure on the category
of modules over an algebra, in such a way that we obtain a
monoidal category. From a purely algebraic point of view, it is
possible to construct quasi-Hopf algebras out of a Hopf algebra
$H$, by twisting the comultiplication: take an invertible element
$F=\sum F^1\ot F^2\in H\ot H$ such that $\sum\va(F^1) F^2=\sum
F^1\va(F^2)=1$, and put $\Delta_F(h)=F\Delta(h)F^{-1}$. Starting
from a finite dimensional cocommutative Hopf algebra and a
Sweedler 3-cocycle $\omega\in H^{\ot 3}$, we can define a
quasi-Hopf algebra structure on the $k$-vector space $H^*\ot H$
\cite{bp}, which is a generalization of the
Dijkgraaf-Pasquier-Roche quasi-Hopf algebra $D^{\omega}(G)$ (where
$H$ is a finite group algebra). An even more general construction
is the quantum double $D(H)$ of a finite dimensional quasi-Hopf
algebra, introduced independently by Majid (\cite{m1}) and Hausser
and Nill (\cite{hn1}, \cite{hn2}). Majid introduced the category
of left Yetter-Drinfeld modules, and used them to describe the
quantum double in the form of an implicit Tannaka-Krein
reconstruction Theorem. Hausser and Nill prove that the category
of $D(H)$-modules is isomorphic to the category of left-right
Yetter-Drinfeld modules.\\%
The Hausser-Nill description of the quantum double uses the
generating matrix formalism. We will reformulate this
intrinsically in \seref{2}; we need this new formulation for our
main result in \seref{3}, and it also has the merit that it
reveals the link between the quantum double construction in the
quasi and classical Hopf algebra case in a transparent way.\\ 
In \seref{3}, we give the quasi-Hopf algebra version of the following
result of Majid \cite{m3}: a finite dimensional Hopf algebra is
quasitriangular if and only if there exists a Hopf algebra
projection of $D(H)$ onto $H$ that covers the natural inclusion.
In this situation, $D(H)$ is a biproduct of a braided Hopf algebra
$B^i$ and $H$, and, as a vector space, $B^i$ is isomophic to $H^*$
(see \cite{m3} for this result in the case of a Hopf algebra). We
give explicit formulas for the structure of $H^*$ as a Hopf
algebra in the category of Yetter-Drinfeld modules. Compare this
to the construction in \cite{bn}, where a braided Hopf algebra
$\un{H}$ is associated to any quasitriangular (not necessarily
finite dimensional) quasi-Hopf algebra $H$, using the Braided
Reconstruction Theorem \cite{m4}.\\ Finally observe that, as in
\cite[Example 3.6]{m3}, for any finite dimensional  quasi-Hopf
algebra $H$, we have an infinite tower of projections
\begin{center}%
\setlength{\unitlength}{1mm}%
\begin{picture}(108,25)(0,0)%
\put(0,20){$H$}%
\put(4, 22){\vector(2,0){10}}%
\put(5,20.5){\vector(-1,0){1}}
\put(5,19.5){$-$}
\put(8,19.5){$-$}
\put(11,19.5){$-$}
\put(15,20){$D(H)$}%
\put(8,23){$i$}%
\put(8,18){$\pi $}%
\put(26,22){\vector(2,0){10}}
\put(36,21){\vector(-2,0){10}}
\put(37,20){$D(D(H))$}
\put(55,22){\vector(2,0){10}}
\put(65,21){\vector(-2,0){10}}
\put(66,20){$D(D(D(H)))$}
\put(91,22){\vector(2,0){10}}
\put(101,21){\vector(-2,0){10}}
\put(103,20){$\cdots $}
\put(30,23){$i_1$}
\put(30,18){$\pi _1$}
\put(59,23){$i_2$}
\put(59,18){$\pi _2$}
\put(95,23){$i_3$}
\put(95,18){$\pi _3$}
\end{picture}
\end{center}%
${\;\;\;}$\\[-25mm]%
The first inclusion is itself covered by a projection if and only if
$H$ is quasitriangular.

\section{Preliminary results}\selabel{1}
\subsubsection*{Quasi-Hopf algebras}
We work over a commutative field $k$. All algebras, linear
spaces etc. will be over $k$; unadorned $\ot $ means $\ot_k$.
Following Drinfeld \cite{d1}, a quasi-bialgebra is
a fourtuple $(H, \Delta
, \va , \Phi )$ where $H$ is an associative algebra with unit,
$\Phi$ is an invertible element in $H\ot H\ot H$, and
$\Delta :\
H\ra H\ot H$ and $\va :\ H\ra k$ are algebra homomorphisms satisfying
the identities
\begin{equation}\label{q1}
(id \ot \Delta )(\Delta (h))=
\Phi (\Delta \ot id)(\Delta (h))\Phi ^{-1},
\end{equation}
\begin{equation}\label{q2}
(id \ot \va )(\Delta (h))=h\ot 1,
\mbox{${\;\;\;}$}
(\va \ot id)(\Delta (h))=1\ot h,
\end{equation}
for all $h\in H$, and $\Phi$ has to be a normalized $3$-cocycle,
in the sense that
\begin{equation}\label{q3}
(1\ot \Phi)(id\ot \Delta \ot id)
(\Phi)(\Phi \ot 1)=
(id\ot id \ot \Delta )(\Phi )
(\Delta \ot id \ot id)(\Phi ),
\end{equation}
\begin{equation}\label{q4}
(id \ot \va \ot id )(\Phi )=1\ot 1\ot 1.
\end{equation}
The map $\Delta $ is called the coproduct or the comultiplication,
$\va $ the counit and $\Phi $ the reassociator. As for Hopf
algebras \cite{sw} we denote $\Delta (h)=\sum h_1\ot h_2$, but
since $\Delta$ is only quasi-coassociative we adopt the further
convention %
$$%
(\Delta \ot id)(\Delta (h))= \sum h_{(1, 1)}\ot h_{(1, 2)}\ot h_2,%
\mbox{${\;\;\;}$}
(id\ot \Delta )(\Delta (h))= \sum h_1\ot h_{(2,
1)}\ot h_{(2,2)}, %
$$%
for all $h\in H$. We will denote the tensor components of $\Phi$
by capital letters, and the ones of $\Phi^{-1}$ by small letters,
namely
\begin{eqnarray*}
&&\Phi=\sum X^1\ot X^2\ot X^3= \sum T^1\ot T^2\ot T^3=
\sum V^1\ot V^2\ot V^3=\cdots\\
&&\Phi^{-1}=\sum
x^1\ot x^2\ot x^3= \sum t^1\ot t^2\ot t^3= \sum v^1\ot v^2\ot v^3=\cdots
\end{eqnarray*}
$H$ is
called a quasi-Hopf algebra if, moreover, there exists an
anti-automorphism $S$ of the algebra $H$ and elements $\a , \b \in
H$ such that, for all $h\in H$, we have:
\begin{equation}\label{q5}
\sum S(h_1)\a h_2=\va (h)\a
\mbox{${\;\;\;}$ and ${\;\;\;}$}
\sum h_1\b S(h_2)=\va (h)\b ,
\end{equation}
\begin{equation} \label{q6}
\sum X^1\b S(X^2)\a X^3=1
\mbox{${\;\;\;}$ and${\;\;\;}$}
\sum S(x^1)\a x^2\b S(x^3)=1.
\end{equation}
For a quasi-Hopf algebra the antipode is determined uniquely up
to a transformation $\a \mapsto U\a $, $\b \mapsto \b U^{-1}$,
$S(h)\mapsto US(h)U^{-1}$, where $U\in H$ is invertible.
The axioms for a quasi-Hopf algebra imply that $\va (\a )\va (\b )=1$,
so, by rescaling $\a $ and $\b $, we may assume without loss of generality
that $\va (\a )=\va (\b )=1$ and $\va \circ S=\va $. The identities
(\ref{q2}), (\ref{q3}) and (\ref{q4}) also imply that
\begin{equation}\label{q7}
(\va \ot id\ot id)(\Phi )=
(id \ot id\ot \va )(\Phi )=1\ot 1\ot 1.
\end{equation}
Next we recall that the definition of a quasi-Hopf
algebra is ``twist coinvariant" in the following sense. An
invertible element $F\in H\ot H$ is called a {\sl gauge transformation}
or {\sl twist} if
$(\va \ot id)(F)=(id\ot \va)(F)=1$.
If $H$ is a quasi-Hopf algebra and $F=\sum F^1\ot F^2\in H\ot H$
is a gauge transformation with inverse $F^{-1}=\sum G^1\ot
G^2$, then we can define a new quasi-Hopf algebra $H_F$ by keeping
the multiplication, unit, counit and antipode of $H$ and
replacing the comultiplication, antipode and the elements $\alpha$
and $\beta$ by
\begin{equation} \label{g1}
\Delta _F(h)=F\Delta (h)F^{-1},
\end{equation}
\begin{equation} \label{g2}
\Phi_F=(1\ot F)(id \ot \Delta )(F) \Phi (\Delta \ot id)
(F^{-1})(F^{-1}\ot 1),
\end{equation}
\begin{equation} \label{g3}
\a_F=\sum S(G^1)\a G^2,
\mbox{${\;\;\;}$}%
\b_F=\sum F^1\b S(F^2).
\end{equation}
It is well-known that the antipode of a Hopf algebra is an
anti-coalgebra morphism. For a quasi-Hopf algebra, we have
the following statement: there exists a gauge transformation
$f\in H\ot H$ such that
\begin{equation} \label{ca}
f\Delta (S(h))f^{-1}= \sum (S\ot S)(\Delta ^{\rm cop}(h))
\mbox{,${\;\;\;}$for all $h\in H$,}
\end{equation}%
where $\Delta ^{\rm cop}(h)=\sum h_2\ot h_1$. $f$ can be computed
explicitly. First set
\begin{equation}
\sum A^1\ot A^2\ot A^3\ot A^4=
(\Phi \ot 1) (\Delta \ot id\ot id)(\Phi ^{-1}),
\end{equation}
\begin{equation} \sum B^1\ot B^2\ot B^3\ot B^4=
(\Delta \ot id\ot id)(\Phi )(\Phi ^{-1}\ot 1)
\end{equation}
and then define $\gamma, \delta\in H\ot H$ by
\begin{equation} \label{gd}%
\gamma =\sum S(A^2)\a A^3\ot S(A^1)\a A^4~~{\rm and}~~
\delta
=\sum B^1\b S(B^4)\ot B^2\b S(B^3).
\end{equation}
$f$ and $f^{-1}$ are then given by the formulas
\begin{eqnarray}
f&=&\sum (S\ot S)(\Delta ^{op}(x^1)) \gamma \Delta (x^2\b
S(x^3)),\label{f}\\ f^{-1}&=&\sum \Delta (S(x^1)\a x^2) \delta
(S\ot S)(\Delta^{\rm cop}(x^3)).\label{g}
\end{eqnarray}
$f$ satisfies the following relations:
\begin{equation} \label{gdf}%
f\Delta (\a )=\gamma ,
\mbox{${\;\;\;}$}
\Delta (\b )f^{-1}=\delta .
\end{equation}
Furthermore the corresponding
twisted reassociator (see (\ref{g2})) is given by
\begin{equation} \label{pf}
\Phi _f=\sum (S\ot S\ot S)(X^3\ot X^2\ot X^1).
\end{equation}
In a Hopf algebra $H$, we obviously have the identity %
$$%
\sum h_1\ot h_2S(h_3)=h\ot 1,~{\rm for~all~}h\in H.%
$$%
We will need the generalization of this formula to the quasi-Hopf
algebra setting. Following \cite{hn1}, \cite{hn2}, we define
\begin{equation} \label{qr}
p_R=\sum p^1\ot p^2=\sum x^1\ot x^2\b S(x^3),%
\mbox{${\;\;\;}$}%
q_R=\sum q^1\ot q^2=\sum X^1\ot S^{-1}(\a X^3)X^2.%
\end{equation}
For all $h\in H$, we then have
\begin{equation} \label{qr1}
\sum \Delta (h_1)p_R[1\ot S(h_2)]=p_R[h\ot 1],
\mbox{${\;\;\;}$}%
\sum [1\ot S^{-1}(h_2)]q_R\Delta (h_1)=(h\ot 1)q_R,%
\end{equation}
and
\begin{equation} \label{pqr}
\sum \Delta (q^1)p_R[1\ot S(q^2)]=1\ot 1,
\mbox{${\;\;\;}$}
\sum [1\ot S^{-1}(p^2)]q_R\Delta (p^1)=1\ot 1,
\end{equation}
\begin{eqnarray}
&&\hspace*{-2cm}(q_R\ot 1)(\Delta \ot id)(q_R)\Phi ^{-1}\nonumber\\
&=&\sum [1\ot S^{-1}(X^3)\ot S^{-1}(X^2)] [1\ot S^{-1}(f^2)\ot
S^{-1}(f^1)] (id \ot \Delta )(q_R\Delta (X^1)),\label{qr2}\\
&&\hspace*{-2cm}\Phi (\Delta \ot id)(p_R)(p_R\ot id)\nonumber\\
&=&\sum (id\ot \Delta )(\Delta (x^1)p_R)(1\ot f^{-1})(1\ot S(x^3)\ot
S(x^2))\label{pr1},
\end{eqnarray}
where $f=\sum f^1\ot f^2$ is the twist defined in (\ref{f}).\\
A quasi-bialgebra or quasi-Hopf algebra $H$ is
quasitriangular if there exists an element
$R=\sum R^1\ot R^2=\sum r^1\ot r^2\in H\ot H$ such that
\begin{eqnarray}
(\Delta \ot id)(R)&=&\sum X^2R^1x^1Y^1\ot X^3x^3r^1Y^2\ot
X^1R^2x^2r^2Y^3\label{qt1}\\ (id \ot \Delta )(R)&=&\sum
x^3R^1X^2r^1y^1\ot x^1X^1r^2y^2\ot x^2R^2X^3y^3\label{qt2}\\
\Delta ^{\rm op}(h)R&=&R\Delta (h),~{\rm for~all~}h\in
H\label{qt3}\\ (\va \ot id)(R)&=&(id\ot \va)(R)=1.\label{qt4}
\end{eqnarray}
In {\cite{bn2}}, it is shown that $R$ is
invertible. Furthermore, the element
\begin{equation} \label{elmu}
u=\sum S(R^2p^2)\a R^1p^1
\end{equation}
(with $p_R=\sum p^1\ot p^2$ defined as in (\ref{qr})) is
invertible in $H$, and
\begin{equation} \label{inelmu}
u^{-1}=\sum X^1R^2p^2S(S(X^2R^1p^1)\a X^3),
\end{equation}
\begin{equation} \label{sqina}
\va (u)=1~~{\rm and}~~
S^2(h)=uhu^{-1}
\end{equation}
for all $h\in H$. Consequently the antipode $S$ is bijective and,
as in the Hopf algebra case, the assumptions about invertibility
of $R$ and bijectivity of $S$ can be deleted. Moreover, the
$R$-matrix $R=\sum R^1\ot R^2$ satisfies the identity (see
\cite{ac}, \cite{hn2}, \cite{bn2}):
\begin{equation} \label{ext}
f_{21}Rf^{-1}=(S\ot S)(R)
\end{equation}
where $f=\sum f^1\ot f^2$ is the twist defined in (\ref{f}), and
$f_{21}=\sum f^2\ot f^1$.\\%
Assume that $(H, \Delta, \va , \Phi)$
is a quasi-bialgebra. For left $H$-modules $U, V, W$, we define a
left $H$-action on $U\ot V$ by %
$$%
h\cdot (u\ot v)=\sum h_1\cd u\ot h_2\cd v.%
$$%
We have isomorphisms $a_{U, V, W}:\ (U\ot V)\ot W\ra
U\ot (V\ot W)$ in ${}_H{\cal M}$ given by
\begin{equation} \label{as}
a_{U, V, W}((u\ot v)\ot w)= \Phi \cd (u\ot (v\ot w)).
\end{equation}
The counit $\varepsilon:\ H\to k$ makes $k\in {}_H{\cal M}$, and
the natural isomorphisms $\lambda:\ k\ot H\to H$ and $\rho:\ H\ot
k\to H$ are in ${}_H{\cal M}$.\\%
With this notation, $({}_H{\cal M}, \ot, k, a, \lambda, \rho)$ is
a monoidal category (see \cite{k} or \cite{Mclane} for the precise
definition).\\ Now let $H$ be a quasitriangular quasi-bialgebra,
with $R$-matrix $R=\sum R^1\ot R^2$. For two left $H$-modules $U$
and $V$, we define $$c_{U,V}:\ U\ot V\to V\ot U$$ by
\begin{equation}\label{br}
c_{U, V}(u\ot v)=\sum R^2\cd v\ot R^1\cd u
\end{equation}
and then $({}_H{\cal M}, \ot, k, a, \lambda, \rho, c)$ is a braided
monoidal category (cf. \cite{k} or \cite{Mclane}).

\subsubsection*{The biproduct for quasi-Hopf algebras}
First we recall from \cite{m1} the notion of Yetter-Drinfeld
module over a quasi-bialgebra.

\begin{definition}
Let $H$ be a quasi-bialgebra with reassociator $\Phi$. A
left $H$-module $M$ together with a left $H$-coaction
$$\lambda_M:\ M\to H\ot M,~~\lambda_M(m)=\sum m_{(-1)}\ot m_{(0)}$$
is called a left Yetter-Drinfeld module if the following
equalities hold, for all $h\in H$ and $m\in M$:
\begin{eqnarray}
&&\hspace*{-2cm} \sum X^1m_{(-1)}\ot (X^2\cd m_{(0)})_{(-1)}X^3
\ot (X^2\cd m_{(0)})_{(0)}\nonumber\\%
&=&\sum X^1(Y^1\cd
m)_{(-1)_1}Y^2\ot X^2(Y^1\cd m)_{(-1)_2}Y^3\ot X^3\cd (Y^1\cd
m)_{(0)}\label{y1}\\%
&&\hspace*{-2cm}\sum \va
(m_{(-1)})m_{(0)}=m\label{y2}\\%
&&\hspace*{-2cm}\sum h_1m_{(-1)}\ot h_2\cd m_{(0)}= \sum (h_1\cd
m)_{(-1)}h_2\ot (h_1\cd m)_{(0)}.\label{y3}
\end{eqnarray}
\end{definition}

The category of left Yetter-Drinfeld $H$-modules and $k$-linear
maps that preserve the $H$-action and $H$-coaction is denoted
$\yd$. In \cite{m1}, it is shown that $\yd$ is a braided monoidal
category. The forgetful functor $\yd\to {}_H{\cal M}$ is monoidal,
and the coaction on the tensor product $M\ot N$ of two
Yetter-Drinfeld modules $M$ and $N$ is given by
\begin{equation}
\lambda _{M\ot N}(m\ot n)=\sum X^1(x^1Y^1\cd m)_
{(-1)} x^2(Y^2\cd n)_{(-1)}Y^3 \ot X^2\cd (x^1Y^1\cd
m)_{(0)}\ot X^3x^3\cd (Y^2\cd n)_{(0)}.\label{y4}
\end{equation}
The braiding is given by
\begin{equation}\label{y5}
c_{M, N}(m\ot n)=\sum m_{(-1)}\cd n\ot m_{(0)}.
\end{equation}
This braiding is invertible if $H$ is a quasi-Hopf algebra \cite{bn},
and its inverse is given by
\begin{equation}
c^{-1}_{M, N}(n\ot m)=\sum y^3_1X^2\cd (x^1\cd
m)_{(0)} \ot S^{-1}(S(y^1)\a y^2X^1(x^1\cd
m)_{(-1)}x^2\b S(y^3_2X^3x^3))\cd n.\label{y6}
\end{equation}
We can consider (co)algebras,
bialgebras and Hopf algebras in the braided category $\yd$.
Let $H$ be a quasi-Hopf algebra, and $B$ a Hopf algebra with bijective
antipode in the category $\yd$. Following \cite{bn}, the $k$-vector
space $B\ot H$ becomes a quasi-Hopf algebra $B\times H$, with
the following structure:
\begin{itemize}
\item as an algebra, $B\times H$ is the smash product from \cite{bpv},
that is, the unit is $1\times 1$, and the multiplication is given by
\begin{equation}\label{sm}
(b\times h)(b^{'}\times h^{'})=\sum (x^1\cd b)(x^2h_1\cd
b^{'})\times x^3h_2h^{'}
\end{equation}
\item the comultiplication, counit and reassociator are given by the
formulas
\begin{eqnarray}
&&\hspace*{-2cm}\Delta (b\times h)=\sum y^1X^1\cd b_1\nonumber\\
&\times &
y^2Y^1(x^1X^2\cd
b_2)_{(-1)}x^2X^3_1h_1\ot y^3_1Y^2\cd (x^1X^2\cd b_2)_{(0)}\times
y^3_2Y^3x^3X^3_2h_2\label{com1}\\
&&\hspace*{-2cm}
\va (b\times h)=\va (b)\va (h)~~;~~
\Phi _{B\times H}=\sum 1\times X^1\ot 1\times
X^2\ot 1\times X^3\label{com2}
\end{eqnarray}
\item the antipode is given by
\begin{equation}\label{ant1}
s(b\times h)=\sum (1\times S(X^1x^1_1b_{(-1)}h)\a
)(X^2x^1_2\cd
S_B(b_{(0)})\times X^3x^2\b S(x^3))
\end{equation}
and
\begin{equation}\label{ant2}
\a _{B\times H}=1\times \a
~~;~~
\b _{B\times H}=1\times \b.
\end{equation}
\end{itemize}
As the reader probably expects, we recover Radford's biproduct \cite{r}
in the situation where
$H$ is a bialgebra or a Hopf algebra.

\section{The quantum double $D(H)$}\selabel{2}
From \cite{hn2}, we recall the definition of the quantum double
$D(H)$ of a finite dimensional quasi-Hopf algebra $H$. Let
$\{e_i\}_{i=\ov {1, n}}$ be a basis of $H$, and $\{e^i\}_{i=\ov
{1, n}}$ the corresponding dual basis of $H^*$. $H^*$ is a
coassociative coalgebra, with comultiplication %
$$%
\widehat {\Delta }(\v )= \sum \v _1\ot \v _2= \sum \limits _{i,
j=1}^n\v (e_ie_j)e^i\ot e^j,%
$$%
or, equivalently, %
$$%
\widehat {\Delta }(\v )=\sum \v _1\ot \v_2 \Leftrightarrow \v
(hh^{'})=\sum \v _1(h)\v _2(h^{'}), %
\mbox{${\;\;\;}$$\forall h, h^{'}\in H$.}%
$$%
$H^*$ is also an $H$-bimodule, by %
$$%
<h\rh \v , h^{'}>=\v (h^{'}h),
\mbox{${\;\;\;}$} %
<\v \lh h, h^{'}>=\v (hh^{'}).%
$$%
We also define, for $\varphi\in H^*$ and $h\in H$: %
$$%
\v\rh h=\sum \v (h_2)h_1, %
\mbox{${\;\;\;}$} %
h\lh \v =\sum \v (h_1)h_2. %
$$%
The convolution is a multiplication on $H^*$; it is not
associative, but only quasi-associative:
$$%
[\v \psi]\xi=\sum (X^1\rh \v \lh x^1)[(X^2\rh \psi \lh x^2)
(X^3\rh \xi \lh x^3)],
\mbox{${\;\;\;}$$\forall \v , \psi , \xi \in H^*$.}%
$$%
We also introduce $\ov {S}:\ H^*\ra H^*$ as the coalgebra
antimorphism dual to $S$, i.e. $<\ov {S}(\v ), h>=<\v , S(h)>$.\\
Now consider
$\Omega \in H^{\ot 5}$ given by
\begin{eqnarray}
&&\hspace*{-2cm} \Omega =\sum \Omega^1\ot \Omega^2\ot \Omega^3\ot
\Omega^4\ot \Omega^5 \nonumber \\%
&=&\sum X^1_{(1, 1)}y^1x^1\ot X^1_{(1, 2)}y^2x^2_1\ot
X^1_2y^3x^2_2\ot \smi (f^1X^2x^3)\ot \smi (f^2X^3), \label{O}
\end{eqnarray}
where $f\in H\ot H$ is the twist defined in (\ref{f}). We define
the quantum double $D(H)=H^*\bowtie H$ as follows: as a $k$-linear
space, $D(H)$ equals $H^*\ot H$, and the multiplication is given
by
\begin{eqnarray}
&&\hspace*{-2cm}(\v \bowtie h)(\psi \bowtie h')\nonumber\\ &=&\sum
[(\Omega ^1\rh \v \lh \Omega ^5)(\Omega ^2\rh \psi _2\lh \Omega
^4)]\bowtie \Omega ^3[(\ov {S}^{-1}(\psi _1)\rh h)\lh \psi
_3]h'\nonumber\\ &=&\sum [(\Omega ^1\rh \v \lh \Omega ^5)(\Omega
^2h_{(1, 1)}\rh \psi \lh \smi (h_2)\Omega ^4)]\bowtie \Omega
^3h_{(1, 2)}h'.\label{mdd}
\end{eqnarray}
It is easy to see that
\begin{equation}\label{mdd1}
(\va \bowtie h)(\v \bowtie h')=\sum h_{(1, 1)}\rh \v \lh \smi (h_2)\bowtie
h_{(1, 2)}h'
\end{equation}
and $(\v \bowtie h)(\va \bowtie h')=\v \bowtie hh'$ for all $\v \in H^*$
and $h, h'\in H$. Following \cite{hn1}, \cite{hn2}, $D(H)$ is an
associative algebra with unit $\va \bowtie 1$, and $H$ is a unital
subalgebra via the morphism $i_D:\ H\ra D(H)$, $i_D(h)=\va \bowtie
h$. In \cite[Theorem 3.9]{hn2}, it is shown that we have a
quasitriangular quasi-Hopf algebra %
$$(D(H), \Delta_D, \va _D, \Phi
_D, S_D, \a _D, \b _D, R_D)$$%
Let us describe the structure; we first introduce
\begin{equation}\label{Ddd}
{\bf D}=\sum \limits _{i=1}^n\smi (p^2)e_ip^1_1\ot
(e^i\bowtie p^1_2)\in H\ot D(H),
\end{equation}
where $p_R=\sum p^1\ot p^2$ id defined as in (\ref{qr}). Then we
introduce the following notation: for $E=\sum \limits _iE_i^1 \ot
E_i^2\ot \cdots \ot E_i^n\in H^{\ot n}$ and $n\leq m$, we let
$E^{n_1n_2\cdots n_n}$ be the element in $H^{\ot m}$ having
$E_i^k$ in the $n_k^{th}$ tensor position, $1\leq k\leq n$, and
$1$ elsewhere.\\%
The comultiplication, counit, reassociator, $R$-matrix and
antipode are then defined for all $h\in H$ by the formulas
\begin{eqnarray}
&&\hspace*{-3cm}\Delta _D(i_D(h))=(i_D\ot i_D)(\Delta
(h))\label{cdd11}\\%
 &&\hspace*{-3cm}(i_D\ot \Delta _D)({\bf
D})=(\Phi _D^{-1})^{231}{\bf D}^{13} \Phi _D^{213}{\bf D}^{12}\Phi
^{-1}_D\label{cdd12}\\%
&&\hspace*{-3cm}\va _D(i_D(h))=\va (h)\label{codd11}\\%
&&\hspace*{-3cm}(id\ot \va _D)({\bf D})=(\va \ot id)({\bf D})\\
&&\hspace*{-3cm}\Phi _D=(i_D\ot i_D\ot i_D)(\Phi)\label{pdd}\\
&&\hspace*{-3cm}R_D=(i_D\ot id)({\bf D})=\sum \limits _{i=1}^n
(\va \bowtie \smi (p^2)e_ip^1_1)\ot (e^i\bowtie p^1_2)\label{rdd}\\
&&\hspace*{-3cm}S_D(i_D(h))=i_D(S(h))\label{add11}\\
&&\hspace*{-3cm}(S\ot S_D)({\bf D})=(id\ot i_D)(f_{21}){\bf
D}(id\ot i_D)(f^{-1}) \label{add12}\\ &&\hspace*{-3cm}\a _D=i_D(\a
), \mbox{${\;\;\;}$} \b _D=i_D(\b )\label{abdd}
\end{eqnarray}
where $f=\sum f^1\ot f^2$ is the twist defined in (\ref{f}) with
inverse $f^{-1}$ and $f_{21}=\sum f^2\ot f^1$. We point out that
the comultiplication is completely determined by (\ref{cdd11}) and
(\ref{cdd12}). Using (\ref{mdd}) and (\ref{cdd12}), we find that
\begin{eqnarray*}
&&\hspace*{-5mm}\sum \limits _{i=1, n}^n(\va \bowtie S^{-1}(p^2)e_ip^1_1)\ot
\Delta _D(e^i\bowtie p^1_2)\\
&=&\sum \limits _{i, j=1}^n [\va \bowtie x^3\smi (p^2)e_ip^1_1
X^2\smi (P^2)e_jP^1_1y^1]\ot [(\va \bowtie x^1X^1)(e^j\bowtie P^1_2y^2)]\ot
[(\va \bowtie x^2)(e^i\bowtie p^1_2X^3y^3)],
\end{eqnarray*}
where $\sum P^1\ot P^2$ is another copy of $p_R$.
Clearly this is equivalent to
\begin{eqnarray}
&&\hspace*{-2cm}
\Delta _D(p^1_1\rh \v \lh \smi(p^2)\bowtie p^1_2)\nonumber\\
&=&\sum
[(\va \bowtie x^1X^1)
(P^1_1y^1\rh \v _2\lh X^2\smi (P^2)
\bowtie P^1_2y^2)]\nonumber\\
&&\hspace*{10mm}\ot
[(\va \bowtie x^2)(p^1_1\rh \v _1\lh x^3\smi (p^2)\bowtie
p^1_2X^3y^3)].\label{cdd13}
\end{eqnarray}
On the other hand, by (\ref{pqr}) we have that
\begin{equation} \label{for1}
\sum (\va \bowtie q^1)(p^1_1\rh \v \lh q^2\smi (p^2)\bowtie
p^1_2)=\v \bowtie 1,
\end{equation}
where $q_R=\sum q^1\ot q^2$ is defined by (\ref{qr});
using (\ref{qr}), (\ref{q3}) and (\ref{q7}) we then obtain the
following relation:
\begin{equation}\label{for2}
\sum q^1_1x^1\ot q^1_2x^2\ot q^2x^3=\sum
Y^1\ot q^1Y^2_1\ot \smi (Y^3)q^2Y^2_2,
\end{equation}
and, using (\ref{for1}), (\ref{cdd11})
and (\ref{cdd13}), we compute that
\begin{eqnarray*}
\Delta _D(\v \bowtie 1)
&=&\sum [(\va \bowtie q^1_1x^1X^1)(P^1_1y^1\rh \v _2\lh
X^2\smi (P^2)\bowtie P^1_2y^2)]\\
&&\hspace*{10mm}\ot [(\va \bowtie q^1_2x^2)(p^1_1\rh
\v _1\lh q^2x^3\smi (p^2)\bowtie p^1_2X^3y^3)]\\
{\rm (\ref{for2})}~~~~&=&\sum [(\va \bowtie Y^1X^1)(P^1_1y^1\rh \v _2\lh
X^2\smi (P^2)\bowtie
P^1_2y^2)]\\
&&\hspace*{10mm}\ot [(\va \bowtie q^1Y^2_1)(p^1_1\rh \v _1\lh \smi (Y^3)q^2
Y^2_2\smi (p^2)\bowtie p^1_2X^3y^3)]\\
{\rm (\ref{mdd1},\ref{qr1})}~~~~&=&\sum [(\va \bowtie Y^1X^1)(P^1_1y^1\rh \v
_2\lh X^2\smi
(P^2)\bowtie P^1_2y^2)]\\
&&\hspace*{10mm}\ot [(\va \bowtie q^1)
(p^1_1Y^2_1\rh \v _1\lh \smi (Y^3)q^2\smi
(p^2)\bowtie p^1_2Y^2_2X^3y^3)]\\
{\rm (\ref{mdd1},\ref{pqr})}~~~~&=&\sum [(\va \bowtie Y^1X^1)(P^1_1y^1\rh \v
_2\lh X^2\smi
(P^2)
\bowtie P^1_2y^2)]\\
&&\hspace*{10mm}\ot (Y^2_1\rh \v _1\lh \smi (Y^3)\bowtie Y^2_2X^3y^3)
\end{eqnarray*}
Combining this with (\ref{cdd11}), we obtain the following explicit formula
for the comultiplication on $D(H)$:
\begin{eqnarray}
&&\hspace*{-2cm}\Delta _D(\v \bowtie h)=
(\va \bowtie X^1Y^1)
(p^1_1x^1\rh \v _2\lh Y^2\smi (p^2)\bowtie p^1_2x^2h_1)\nonumber\\
&&\hspace*{10mm}\ot (X^2_1\rh \v _1\lh \smi (X^3)\bowtie X^2_2Y^3x^3h_2)
\label{cddf}
\end{eqnarray}
Similar computations yield explicit formulas for the counit and the
antipode, namely
\begin{equation}\label{coddf}
\va _D(\v \bowtie h)=\va (h)\v (\smi (\a ))
\end{equation}
and
\begin{equation}\label{anddf}
S_D(\v \bowtie h)=\sum (\va \bowtie S(h)f^1)(p^1_1U^1\rh \ov
{S}^{-1}(\v )\lh f^2\smi (p^2)\bowtie p^1_2U^2).
\end{equation}
Here the twist $f=\sum f^1\ot f^2$, its inverse
$f^{-1}=\sum g^1\ot g^2$ and $q_R=\sum q^1\ot q^2$
are defined as in (\ref{f}), (\ref{g}) and (\ref{qr}),
and
\begin{equation}\label{U}
U=\sum U^1\ot U^2=\sum g^1S(q^2)\ot g^2S(q^1).
\end{equation}

\section{$D(H)$ when $H$ is quasitriangular}\selabel{3}
In \cite{m3} it is shown that a finite dimensional Hopf algebra is
quasistriangular if and only if its Drinfeld double is a Hopf
algebra with a projection. In this Section, we will generalize
this result to quasi-Hopf algebras. In this case, $D(H)$ is
isomorphic to a biproduct in the sense of \cite{bn}, between a
certain Hopf algebra $B^i$ in the braided category $\yd $ and $H$.
We will show that, as in the case of a classical Hopf algebra,
$B^i$ equals $H^*$ as a vector space, but with a different
multiplication and comultiplication: the structures of $H^*$ in
$\yd $ are induced by the $R$-matrix and the quasi-Hopf algebra
structure of $H$.\\%
Let $(H,\Delta , \Phi, \va , S, \a , \b )$ be
a quasi-Hopf algebra and $(A, \Delta _A, \Phi _A)$ a
quasi-bialgebra. Recall that a quasi-bialgebra map $\nu:\ H\to A$
is a $k$-algebra map that preserves comultiplication and counit,
such that $\Phi_A=\nu^{\ot 3}(\Phi )$. If $(A, \Delta _A, \Phi
_A,S_A, \a_A , \b_A )$ is also a quasi-Hopf algebra, then a
quasi-bialgebra map $\nu$ is a Hopf algebra map if $\nu (\a )=\a
_A$, $\nu (\b )=\b _A$ and $S_A\circ \nu =\nu \circ S$. Observe
that a bialgebra map between classical Hopf algebras is
automatically a Hopf algebra map. If $\nu$ is a quasi-Hopf algebra
map, then the elements $\gamma,\delta,f,f^{-1}\in H\ot H$ are
mapped by $\nu\ot \nu$ to the corresponding ones in $A\ot A$.\\
Our first result is a generalization of \cite[Corollary 3.2]{m3}.

\begin{lemma}\lelabel{3.1}
Let $H$ be a finite dimensional quasi-Hopf algebra. Then
there exists a quasi-Hopf algebra
projection $\pi :\ D(H)\ra H$ covering the canonical inclusion
$i_D: H\ra D(H)$ if and only if $H$ is quasitriangular.
\end{lemma}

\begin{proof}
First assume that there is a quasi-Hopf algebra morphism
$\pi :\ D(H)\ra H$ such that $\pi \circ i_D=id_H$. Then it is not
hard to see that $R=\sum \pi (R^1_D)\ot \pi (R^2_D)$ is an
$R$-matrix for $H$, where $R_D=\sum R^1_D\ot R^2_D$ is the
canonical $R$-matrix of $D(H)$ defined in (\ref{rdd}). Thus $H$ is
quasitriangular.\\
Conversely, let $H$ be a quasitriangular quasi-Hopf algebra with
$R$-matrix $R=\sum R^1\ot R^2$, and define $\pi :\ D(H)\ra H$ by
\begin{equation} \label{pdd2}
\pi (\v \bowtie h)=\sum \v (q^2R^1)q^1R^2h,
\end{equation}
where $q_R=\sum q^1\ot q^2$ is as in
(\ref{qr}). We have to show that $\pi $ is a quasi-Hopf algebra
morphism and $\pi \circ i_D=id_H$.\\
As before, we write
$q_R=\sum Q^1\ot Q^2$ and $R=\sum r^1\ot r^2$, and then compute
for all $\v ,\psi \in H^*$ and $h, h'\in H$ that
\begin{eqnarray*}
&&\hspace*{-2cm}\pi ((\v \bowtie h)(\psi \bowtie h'))\\%
{\rm (\ref{mdd})}~~~~ %
&=&\sum \pi ((\O ^1\rh \v \lh \O ^5)(\O ^2h_{(1, 1)}\rh \psi \lh
\smi (h_2)\O ^4)\bowtie \O ^3h_{(1, 2)}h')\\%
&=&\sum \v (\O ^5q^2_1R^1_1\O ^1) \psi (\smi (h_2)\O
^4q^2_2R^1_2\O ^2h_{(1, 1)}) q^1R^2\O ^3h_{(1, 2)}h'\\%
{\rm (\ref{O})}~~~~%
&=&\sum \v (\smi (f^2X^3)q^2_1R^1_1X^1_{(1, 1)}y^1x^1) \psi
(\smi (f^1X^2x^3h_2)q^2_2R^1_2X^1_{(1, 2)}%
y^2x^2_1h_{(1, 1)})\\%
&&\hspace*{1cm}q^1R^2X^1_2y^3x^2_2h_{(1, 2)}h'\\%
{\rm (\ref{qt3},\ref{qr2})}~~~~ %
&=&\sum \v (q^2Q^1_2z^2R^1_1y^1x^1)\psi (\smi
(x^3h_2)Q^2z^3R^1_2y^2x^2_1h_{(1, 1)})\\
&&\hspace*{1cm}q^1Q^1_1z^1R^2y^3x^2_2h_{(1, 2)}h'\\%
{\rm (\ref{qt1})}~~~~%
&=&\sum \v (q^2Q^1_2R^1y^1x^1)\psi (\smi
(x^3h_2)Q^2y^3r^1x^2_1h_{(1, 1)}) q^1Q^1_1R^2y^2r^2x^2_2h_{(1,
2)}h'\\%
{\rm (\ref{qr},\ref{qt3})}~~~~%
&=&\sum \v (q^2R^1X^1_1y^1x^1)\psi (\smi (\a
X^3x^3h_2)X^2y^3x^2_2h_{(1, 2)}r^1)q^1R^2X^1_2y^2x^2_1h_{(1,
1)}r^2h'\\%
{\rm (\ref{q3},\ref{q7})}~~~~%
&=&\sum \v (q^2R^1)\psi (\smi (h_2)Q^2h_{(1, 2)}r^1)
q^1R^2Q^1h_{(1, 1)}r^2h' \\%
{\rm (\ref{qr1})}~~~~%
&=&\pi (\v \bowtie h)\pi (\psi \bowtie h').
\end{eqnarray*}
From (\ref{qt4}) and (\ref{qr}), it follows that $\pi (\va \bowtie
h)=h$, for any $h\in H$. Thus we have shown that $\pi $ is an
algebra map, and that $\pi \circ i_D=id_H$.\\ Since $\pi \circ
i_D=id _H$, we have that $(\pi \ot \pi \ot \pi )(\Phi _D)=\Phi $.
$\pi$ also preserves the comultiplication, since
\begin{eqnarray*}
&&\hspace*{-2cm}\sum \pi ((\v \bowtie h)_1)\ot \pi ((\v \bowtie h)_2)\\
{\rm (\ref{cddf})}~~~~&=&
\sum \pi ((\va \bowtie X^1Y^1)(p^1_1x^1\rh \v _2\lh Y^2\smi
(p^2)\bowtie p^1_2x^2h_1))\\
&&\hspace*{1cm}\ot \pi (X^2_1\rh \v _1\lh \smi (X^3)\bowtie X^2_2Y^3x^3h_2)\\
&=&\sum \v _2(Y^2\smi (p^2)q^2R^1p^1_1x^1) \v _1(\smi
(X^3)Q^2r^1X^2_1)X^1Y^1q^1R^2p^1_2x^2h_1\ot Q^1r^2X^2_2Y^3x^3h_2\\
&&~~~~({\rm since}~\pi~\mbox{is an algebra map and}~\pi \circ i_D=id_H)\\
{\rm (\ref{qt3})}~~~~&=&
\sum \v (\smi (X^3)Q^2r^1X^2_1Y^2\smi
(p^2)q^2p^1_2R^1x^1) X^1Y^1q^1p^1_1R^2x^2h_1\ot Q^1r^2X^2_2Y^3x^3h_2\\
{\rm (\ref{pqr},\ref{qt3})}~~~~&=&
\sum \v (\smi (X^3)Q^2X^2_2r^1Y^2R^1x^1)X^1Y^1R^2x^2h_1
\ot Q^1X^2_1r^2Y^3x^3h_2\\
{\rm (\ref{for2})}~~~~&=&
\sum \v (q^2y^3r^1Y^2R^1x^1)q^1_1y^1Y^1R^2x^2h_1\ot
q^1_2y^2r^2Y^3x^3h_2\\
{\rm (\ref{qt2})}~~~~&=&\sum \v (q^2R^1)\Delta (q^1R^2h)=%
\Delta (\pi (\v \bowtie h)).
\end{eqnarray*}
$\pi$ also preserves the counit, since $(\va \circ \pi )(\v \bowtie h) =\v
(\smi (\a ))\va
(h)=\va _D(\v \bowtie h)$.\\
We easily see that $\pi
(\a _D)=\a $ and $\pi (\b _D)=\b $, so we are done if we can show that
$S\circ \pi =\pi \circ S_D$. This follows from the next computation:
\begin{eqnarray*}
&&\hspace*{-2cm}(\pi \circ S_D)(\v \bowtie h)\\
{\rm (\ref{anddf})}~~~~&=&\sum \pi ((\va \bowtie S(h)f^1)(p^1_1U^1\rh \ov
{S}^{-1}
(\v )\lh f^2\smi (p^2)\bowtie p^1_2U^2))\\
&=&\sum <\ov {S}^{-1}(\v ), f^2\smi (p^2)q^2R^1p^1_1U^1>
S(h)f^1q^1R^2p^1_2U^2\\
&&~~~~({\rm since}~\pi~\mbox{is an algebra map and}~\pi \circ i_D=id_H)\\
{\rm (\ref{qt3},\ref{pqr},\ref{U})}~~~~&=&\sum <\ov {S}^{-1}(\v ),
f^2R^1g^1S(q^2)>S(h)f^1R^2g^2S(q^1)\\
{\rm (\ref{ext})}~~~~&=&
\sum <\ov {S}^{-1}(\v ), S(q^2R^1)>S(q^1R^2h)\\
&=&(S\circ \pi)(\v \bowtie h),
\end{eqnarray*}
and this completes our proof.
\end{proof}

Observe that, in the case where $(H, R)$ is quasitriangular, the map
$\pi $ given by (\ref{pdd2}) is a quasitriangular morphism, i.e.
$(\pi \ot \pi )(R_D)=R$. Indeed,
\begin{eqnarray*}
(\pi \ot \pi )(R_D)
&=&\sum \limits _{i=1}^n<\pi \ot \pi , \va \bowtie
\smi (p^2)e_ip^1_1\ot e^i\bowtie p^1_2>
\mbox{${\;\;\;}$by (\ref{rdd})}\\
{\rm (\ref{pdd2})}~~~~&=&\sum \limits _{i=1}^n
e^i(q^2R^1)\smi (p^2)e_ip^1_1\ot q^1R^2p^1_2\\
{\rm (\ref{qt3})}~~~~&=&\sum \smi (p^2)q^2p^1_2R^1\ot q^1p^1_1R^2\\%
{\rm (\ref{pqr})}~~~~&=&R.
\end{eqnarray*}

The structure of a quasi-Hopf algebra with projection
was given in \cite{bn}. More precisely, if $H$ is a
quasi-Hopf algebra, and $A$ a quasi-bialgebra with two
quasi-bialgebra maps
\begin{center}
\setlength{\unitlength}{1mm}
\begin{picture}(35,25)(0,0)
\put(0,20){$H$}
\put(4, 22){\vector(2,0){10}}
\put(14,21){\vector(-2,0){10}}
\put(15,20){$A$}
\put(8,23){$i$}
\put(8,18){$\pi $}
\put(79,20){$(3.2)$}
\end{picture}
\end{center}
${\;\;}$\\[-25mm]
such that $\pi \circ i=id$, then there exists a
braided bialgebra $B^i$ in the category $\yd $ such that $A\cong
B^i\times H$ as quasi-bialgebras. In fact, if we define on $A$ a
new multiplication given by
\begin{equation} \label{ima}
a\circ a'=\sum i(X^1)ai(S(x^1X^2)\a
x^2X^3_1)a'i(S(x^3X^3_2))
\end{equation}
and we denote this new structure on $A$ by $A^i$, then $A^i$
becomes a left $H$-module algebra (i.e. an algebra in
the monoidal category ${}_H{\cal M}$), with unit $i(\b )$ and with the
left adjoint action induced by $i$, that is $h\tr _ia=\sum
i(h_1)ai(S(h_2))$, for all $h\in H$ and $a\in A$. Now let
\begin{equation} \label{b}
B^i=\{a\in A\mid \sum a_1\ot \pi (a_2)=\sum
i(x^1)ai(S(x^3_2X^3)f^1)\ot x^2X^1\b S(x^3_1X^2)f^2\}.
\end{equation}
First, $B^i$ is a subalgebra of $A^i$ in the monoidal category
$_H{\cal M}$. Secondly, $B^i$ is an algebra in the braided
category $\yd $ where the left coaction of $H$ on $B^i$ is given
by
\begin{equation} \label{lbi}
\l _{B^i}(b)=\sum X^1Y^1_1\pi (b_1)g^1S(q^2Y^2_2)Y^3\ot
i(X^2Y^1_2)b_2i(g^2S(X^3q^1Y^2_1)).
\end{equation}
Here $f^{-1}=\sum g^1\ot g^2$ and $q_R=\sum q^1\ot q^2$ are
as in (\ref{g}) and (\ref{qr}).\\
${\;\;\;}$%
Also, as $k$-vector space, $B^i$ is the image of the
$k$-linear map $\Pi :\ A\ra A$ defined by
\begin{equation} \label{Pi}
\Pi (a)=\sum a_1i(\b S(\pi (a_2))).
\end{equation}
For all $a\in A$, we define
\begin{equation} \label{dbi}
\un {\Delta }(\Pi (a))=\sum \Pi (a_1)\ot \Pi (a_2).
\end{equation}
This makes $B^i$ into a coalgebra in $\yd $ and a
bialgebra in $\yd $. The counit of $B^i$ is $\un {\va }=
\va \mid _B$. The bialgebra isomorphism $\chi :\ B^i\times
H\ra A$ is given by
\begin{equation}\label{th}
\chi (b\times h)=\sum i(X^1)bi(S(X^2)\a X^3h).
\end{equation}
In the situation where $A$ is a quasi-Hopf algebra, with antipode $S_A$,
and $i$
and $\pi $ are quasi-Hopf algebra maps, we have that $B^i$ is a braided
Hopf algebra
in $\yd $ with antipode $\un {S}$ given by
\begin{equation} \label{ant}
\un {S}(b)=\sum i(\pi (b_1)\b )S_A(b_2).
\end{equation}
Moreover, in this case $\un {S}$ is bijective and $\chi$ is a
quasi-Hopf algebra isomorphism.\\
Now we will apply the above description of a quasi-Hopf algebra
with a projection in our case, namely $(H, R)$ is a finite dimensional
quasitriangular quasi-Hopf algebra, $A=D(H)$, $i=i_D$ is the
canonical inclusion and $\pi $ is the map defined by (\ref{pdd2}). As
in the Hopf case, we show first that $B^i$ is isomorphic to $H^*$
as $k$-vector spaces. This follows from
\begin{eqnarray*}
\Pi (\v \bowtie h)
&=&\sum (\v \bowtie h)_1(\va \bowtie \b S(\pi ((\v \bowtie h)_2)))\\
{(\rm \ref{cddf},\ref{pdd2})}~~~~
&=&\sum \v _1(\smi
(X^3)q^2R^1X^2_1)(\va \bowtie X^1Y^1)\\
&&\hspace*{1cm}(p^1_1y^1\rh\v _2\lh Y^2\smi (p^2)\bowtie
p^1_2y^2h_1)
(\va \bowtie \b S(q^1R^2X^2_2Y^3y^3h_2))\\%
{(\rm \ref{mdd},\ref{q5})}~~~~
&=&\va (h)\sum \v _1(\smi
(X^3)q^2R^1X^2_1)(\va \bowtie X^1Y^1)\\%
&&\hspace*{1cm}(p^1_1y^1\rh \v _2\lh Y^2\smi (p^2)\bowtie p^1_2y^2\b
S(q^1R^2X^2_2Y^3y^3))\\%
&=&\va (h)\Pi (\v \bowtie 1)
\end{eqnarray*}
for all $\v \in H^*$ and $h\in H$, so $B^i=\Pi (D(H))=\Pi
(H^*\bowtie 1)$; this means that we have a surjective $k$-linear map
$$%
H^*\to B^i:\ \v\mapsto \Pi(\v\bowtie 1).%
$$%
This map is also injective since
\begin{eqnarray*}
&&\hspace*{-2cm}<id\ot \va , \Pi (\v \bowtie 1)>\\
&=&\sum \v _1(\smi (X^3)q^2R^1X^2_1)\\
&&\hspace*{1cm}<id\ot \va , (\va \bowtie
X^1Y^1)(p^1_1y^1\rh \v _2\lh Y^2\smi (p^2)\bowtie p^1_2y^2\b
S(q^1R^2X^2_2Y^3y^3)>\\
&=&\sum <id\ot \va , (\va \bowtie X^1Y^1)(p^1_1y^1\rh \v \lh \smi
(X^3)q^2R^1X^2_1Y^2\smi (p^2)\\
&&\hspace*{1cm}\bowtie p^1_2y^2\b S(q^1R^2X^2_2Y^3y^3))>\\
{\rm (\ref{mdd})}~~~~&=&\sum <id\ot \va , (X^1Y^1)_{(1, 1)}p^1_1y^1\rh \v
\lh \smi
(X^3)q^2R^1X^2_1Y^2\smi ((X^1Y^1)_2p^2)\\
&&\hspace*{1cm}\bowtie (X^1Y^1)_{(1, 2)}p^1_2y^2\b S(q^1R^2X^2_2Y^3y^3)>\\
{\rm (\ref{qt4})}~~~~&=&\sum X^1_1p^1\rh \v \lh \smi (\a X^3)X^2\smi
(X^1_2p^2)\\%
{\rm (\ref{qr},\ref{pqr})}~~~~&=&\v .
\end{eqnarray*}
In fact we have shown that the map
$$\mu:\ B^i=\Pi (D(H))\ra H^*,~~~\mu (\Pi (\v \bowtie h))=\va (h)\v $$
is an isomorphism of $k$-vector spaces, with inverse
$$%
\mu ^{-1}(\v )=\Pi (\v\bowtie 1).
$$
From now on, $\un {H}^*$ will be the $k$-vector space $H^*$, with
the structure of Hopf algebra in the braided category $\yd $
induced from $B^i$ via $\mu$. Let us compute the structure maps
of $\un {H}^*$ in $\yd $.

\begin{proposition}\prlabel{3.2}
The structure of $\un{H}^*$ as a Yetter-Drinfeld module is given
by the formulas
\begin{eqnarray}
h\cd \v &=&\sum h_1\rh \v \lh \smi (h_2)\label{mhs}\\ \l _{\un
{H}^*}(\v )&=&\sum R^2\ot R^1\cd \v .\label{chs}
\end{eqnarray}
\end{proposition}

\begin{proof}
(\ref{mhs}) is easy, and left to the reader.
Observe that (\ref{q3}) and (\ref{q5}) imply
\begin{equation}\label{Fo2}
\sum X^1p^1_1\ot X^2p^1_2\ot X^3p^2=\sum x^1\ot x^2_1p^1\ot
x^2_2p^2S(x^3)
\end{equation}
and this allows us to compute
\begin{eqnarray}
\Pi (\v \ot 1)
&=& \sum (\va \bowtie X^1Y^1)\nonumber\\
&&\hspace*{1cm}(p^1_1y^1\rh \v \lh \smi (X^3)q^2R^1X^2_1Y^2
\smi (p^2)\bowtie p^1_2y^2\b S(q^1R^2X^2_2Y^3y^3))\nonumber\\
{\rm (\ref{qt3},\ref{for2},\ref{qr})}~~~~&=&
\sum (\va \bowtie q^1_1x^1Y^1)(p^1_1P^1\rh \v \lh q^2x^3R^1Y^2\smi (p^2)
\bowtie p^1_2P^2S(q^1_2x^2R^2Y^3))\nonumber\\
{\rm (\ref{mdd},\ref{pr1})}~~~~&=&
\sum (\va \bowtie q^1_1x^1)(y^1p^1\rh \v \lh q^2x^3R^1\smi (y^3p^2_2g^2)
\bowtie y^2p^2_1g^1S(q^1_2x^2R^2))\nonumber\\
{\rm (\ref{ext},\ref{qt3},\ref{mdd})}~~~~
&=&\sum (q^1_1x^1)_{(1, 1)}y^1p^1\rh \v \lh q^2x^3
\smi ((q^1_1x^1)_2y^3R^1p^2_1g^1)\bowtie (q^1_1x^1)_{(1, 2)}
y^2\nonumber\\
&&\hspace*{1cm}R^2p^2_2g^2S(q^1_2x^2)\nonumber\\
{\rm (\ref{q1},\ref{qt3},\ref{pr1})}~~~~
&=&\sum y^1(q^1_1)_1X^1p^1_1P^1\rh \v \lh q^2\smi
(y^3R^1(q^1_1)_{(2, 1)}X^2p^1_2P^2) \nonumber\\
&&\hspace*{1cm}\bowtie y^2R^2(q^1_1)_{(2, 2)}X^3p^2S(q^1_2)
\nonumber\\
{\rm (\ref{q1},\ref{qr1})}~~~~&=&\sum y^1X^1p^1_1q^1_1P^1\rh \v \lh \smi
(y^3R^1X^2p^1_2q^1_2P^2S(q^2))\bowtie  y^2R^2X^3p^2\nonumber\\
{\rm (\ref{pqr},\ref{Fo2})}~~~~&=&\sum y^1x^1\rh \v \lh \smi %
(y^3R^1x^2_1p^1)\bowtie y^2R^2x^2_2p^2S(x^3). \label{fPv}
\end{eqnarray}
A similar computation, using (\ref{fPv}) yields
\begin{eqnarray}
&&\hspace*{-2cm}
\pi (\Pi (\v \bowtie 1)_1)\ot \Pi (\v \bowtie 1)_2
=\sum X^1Y^1r^2z^2y^2_1R^2_1x^2_{(2,
1)}p^2_1S(x^3)_1\nonumber\\ &\ot& X^2_1Y^2r^1z^1y^1x^1\rh \v \lh
\smi (X^3y^3R^1x^2_1p^1)\bowtie X^2_2Y^3 z^3y^2_2R^2_2x^2_{(2,
2)}p^2_2S(x^3)_2. \label{Fo1}
\end{eqnarray}
Here $R=\sum r^1\ot r^2$ is another copy of $R$. Since
\begin{eqnarray}
&&\hspace*{-2cm}<id\ot \va , \Pi ((\v \bowtie h)(\Psi \bowtie h'))>=
<id \ot \va , (\v \bowtie h)(\Psi \bowtie h')>\\
&=&
\va (h')\sum (X^1_1x^1\rh \v \lh \smi (f^2X^3))(X^1_2x^2h_1\rh \Psi
\lh \smi (f^1X^2x^3h_2))\label{Fo3},
\end{eqnarray}
for all $\v , \Psi \in H^*$ and $h, h'\in H$, we obtain, using
(\ref{lbi}) that
\begin{eqnarray*}
\l _{\un {H}^*}(\v )
&=&\sum q^1_1X^1r^2y^2R^2x^2_2p^2S(q^2x^3) \\
&&\hspace*{1cm}q^1_{(2, 1)}X^2r^1
y^1x^1\rh \v \lh \smi (q^1_{(2, 2)}X^3y^3R^1x^2_1p^1)\\
{\rm (\ref{q1},\ref{qt3})}~~~~&=&\sum X^1r^2q^1_{(1,
2)}y^2R^2x^2_2p^2S(q^2x^3) \\
&&\hspace*{1cm}\ot X^2r^1q^1_{(1, 1)}y^1x^1\rh \v \lh \smi
(X^3q^1_2y^3R^1x^2_1p^1)\\
{\rm (\ref{q1},\ref{qt3},\ref{for2})}~~~~
&=&\sum X^1r^2y^2R^2q^1_2Y^2_{(1, 2)}p^2S(q^2Y^2_2)Y^3\\
&&\hspace*{1cm}\ot
X^2r^1y^1Y^1\rh \v \lh
\lh \smi (X^3y^3R^1q^1_1Y^2_{(1, 1)}p^1)\\
{\rm (\ref{qr1},\ref{pqr})}~~~~
&=&\sum X^1r^2y^2R^2Y^3\ot X^2r^1y^1Y^1\rh \v \lh \smi (X^3y^3R^1Y^2)\\
{\rm (\ref{qt1},\ref{mhs})}~~~~%
&=&\sum R^2\ot R^1\cd \v .
\end{eqnarray*}
\end{proof}

Our next goal is to compute the algebra structure of $\un {H}^*$.
First we need the following Lemma.

\begin{lemma}\lelabel{3.3}
Let $H$ be a quasi-Hopf algebra, $A$ a quasi-bialgebra and
assume that $\pi:\ A\to H$ is a quasi-bialgebra map that is a
left inverse of the quasi-bialgebra map $i:\ H\to A$. Then the map
$\Pi$ defined in (\ref{Pi}) satisfies the formula
\begin{equation} \label{rPi}
\Pi (\Pi (a)\circ a')=\Pi (a)\circ \Pi (a').
\end{equation}
\end{lemma}

\begin{proof}
Recall from \cite{bn} that $\Pi (ai(h))=\va (h)\Pi (a)$, $\Pi
(i(h)a)=h\tr _i\Pi (a)$ and $\Pi (aa')=\sum a_1\Pi
(a')i(S(\pi (a_2)))$, for all $a, a'\in A$, $h\in H$ (see
\cite{bn}). Then we compute
\begin{eqnarray*}
\Pi (\Pi (a)\circ a') &=&\sum \Pi (i(X^1)\Pi
(a)i(S(x^1X^2)\a x^2X^3_1)a'i(S(x^3X^3_2)))\\ {\rm
(\ref{qr})}~~~~&=&\sum q^1\tr _i\Pi (\Pi (a)i(S(q^2))a')\\ &=&\sum
q^1\tr _i[\Pi (a)_1\Pi(i(S(q^2))a')i(S(\pi (\Pi(a)_2)))]\\ &=&\sum
q^1\tr _i[i(x^1)\Pi (a)i(S(x^3_2X^3)f^1)\Pi
(i(S(q^2))a')i(S(x^2X^1\b S(x^3_1X^2)f^2))]\\%
&&~~~~~\mbox{since $\Pi (a)\in B^i$}\\%
{\rm (\ref{ca})}~~~~&=&\sum q^1\tr
_i[i(x^1)\Pi (a)i(S(q^2_2x^3_2X^3)f^1)\Pi (a')i(S(x^2X^1\b
S(q^2_1x^3_1X^2)f^2))]\\ {\rm (\ref{q3},\ref{q5})}~~~~&=&\sum
q^1X^1\tr _i[i(x^1)\Pi (a)i(S(q^2_2X^3)f^1)\Pi
(a')i(S(x^2\b S(q^2_1X^2x^3)f^2))]\\%
{\rm (\ref{qr2},\ref{qr})}~~~~&=&\sum q^1Q^1_1y^1_{(1, 1)}\tr _i[i(p^1)\Pi
(a)i(S(Q^2y^1_2)y^2)\Pi (a')i(S(p^2S(q^2Q^1_2y^1_{(1,
2)})y^3))]\\
{\rm (\ref{qr1})}~~~~&=&\sum i(q^1_1p^1Q^1y^1_1)\Pi (a)i(S(Q^2y^1_2)y^2)\Pi
(a')i(S(q^1_2p^2S(q^2)y^3))\\
{\rm (\ref{pqr})}~~~~&=&\sum i(Q^1y^1_1)\Pi (a)i(S(Q^2y^1_2)y^2)\Pi
(a')i(S(y^3))\\
{\rm (\ref{qr},\ref{q3},\ref{q5})}~~~~&=&\Pi (a)\circ \Pi (a')
\end{eqnarray*}
\end{proof}

\begin{proposition}\prlabel{3.4}
The structure of $\un {H}^*$ as a Hopf algebra in $\yd$ is given by the
formulas
\begin{eqnarray}
&&\hspace*{-2cm}
\v \circ \Psi =\sum (x^1X^1\rh \v \lh \smi (f^2x^3_2Y^3R^1X^2))\nonumber\\
&&\hspace*{1cm}
(x^2Y^1R^2_1X^3_1\rh \Psi \lh \smi (f^1x^3_1Y^2R^2_2X^3_2))\label{mbdd}\\
&&\hspace*{-2cm}
\un {\Delta }(\v )=\sum X^1_1p^1\rh \v _2\lh \smi (X^1_2p^2)\ot X^2
\rh \v _1\lh \smi (X^3)\label{cbdd}\\
&&\hspace*{-2cm}
\un {\va }(\v )=\v (\smi (\a ))\label{cbdd1}\\
&&\hspace*{-2cm}
\un {S}(\v )=\sum Q^1q^1R^2x^2\cd [p^1P^2S(Q^2)\rh
\ov {S}^{-1}(\v )\lh S(q^2R^1x^1P^1)x^3\smi (p^2)]\label{anbdd}
\end{eqnarray}
The unit element is $\va$.
\end{proposition}

\begin{proof}
It follows from (\ref{rPi}) that the multiplication $\circ $ on $\un {H}^*$ is
\begin{eqnarray*}
\v \circ \Psi
&=&<id\ot \va , \Pi (\v \bowtie 1)\circ \Pi (\Psi \bowtie 1)>\\
&=&<id \ot \va , \Pi (\Pi (\v \bowtie 1)\circ (\Psi \circ 1))>.
\end{eqnarray*}
Now
\begin{eqnarray*}
\Pi (\v \bowtie 1)\circ (\Psi \bowtie 1)
&=&\sum i(X^1)\Pi (\v \bowtie 1)i(S(z^1X^2)\a z^2X^3_1)
(\Psi \bowtie 1)i(S(z^3X^3_2))\\
{\rm (\ref{q3},\ref{q5},\ref{fPv})}~~~~&=&\sum i(q^1z^1_1)(y^1x^1\rh \v \lh
\smi
(y^3R^1x^2_1p^1)\bowtie  y^2R^2x^2_2p^2S(x^3))\\
&&\hspace*{1cm} i(S(q^2z^1_2)z^2)(\Psi \bowtie S(z^3))\\
{\rm (\ref{mdd},\ref{q1},\ref{qt3})}~~~~
&=&\sum (y^1q^1_1z^1_{(1, 1)}x^1\rh \v \lh \smi
(y^3R^1q^1_{(2, 1)}(z^1_{(1, 2)}x^2)_1p^1)\bowtie y^2R^2q^1_{(2, 2)}\\
&&\hspace*{1cm}(z^1_{(1, 2)}x^2)_2
p^2S(q^2z^1_2x^3)z^2)(\Psi \bowtie S(z^3))\\
{\rm (\ref{q1},\ref{qr1})}~~~~
&=&\sum [y^1q^1_1x^1z^1_1\rh \v \lh \smi (y^3R^1q^1_{(2, 1)}x^2_1p^1z^1_2)
\bowtie y^2R^2q^1_{(2, 2)}x^2_2p^2\\
&&\hspace*{1cm}S(q^2x^3)z^2](\Psi \bowtie S(z^3))\\
{\rm (\ref{for2})}~~~~
&=&\sum [y^1X^1z^1_1\rh \v \lh \smi (y^3R^1q^1_1X^2_{(1, 1)}p^1z^1_2)\bowtie
y^2R^2q^1_2X^2_{(1, 2)}p^2\\
&&\hspace*{1cm}S(q^2X^2_2)X^3z^2](\Psi \bowtie S(z^3))\\
{\rm (\ref{qr1},\ref{pqr})}~~~~&=&\sum y^1X^1z^1_1\rh \v \lh \smi
(y^3R^1X^2z^1_2)\bowtie %
y^2R^2X^3z^2) (\Psi \bowtie S(z^3)).
\end{eqnarray*}
Using (\ref{Fo3})
and (\ref{q3}) we then obtain (\ref{mbdd}).
It is easy to prove that $\va$ is the unit element of $\un {H}^*$.\\
(\ref{cbdd}-\ref{cbdd1}) follow from the following formula for the
comultiplication
in $D(H)$ (see \cite{bn}), and the defining axioms of a quasi-bialgebra.
We leave further detail to the reader:
$$%
\un {\Delta }(\Pi (\v \bowtie 1))=\sum \Pi (i(X^1)(\v \bowtie 1)_1)\ot
i(X^2)(\v \bowtie 1)_2i(S(X^3)).
$$%
(\ref{anbdd}) follows after a straigthforward, but long and tedious computation
using (\ref{Fo1}), (\ref{anddf}) and (\ref{Fo3}).
\end{proof}

\begin{remarks}\rm
1) The quasi-Hopf algebra isomorphism
$\chi :\ \un {H}^*\times H\ra D(H)$ is given by
\begin{equation} \label{chidd}
\chi (\v \times h)=\sum x^1X^1\rh \v \lh \smi (x^3R^1X^2)\bowtie
x^2R^2X^3h.
\end{equation}

2)
Let $(H, R)$ be a quasitriangular
quasi-Hopf algebra. Then we have a monoidal functor
${\bf F}:\ {}_H{\cal M}\ra \yd $ which sends algebras, coalgebras,
bialgebras etc.
in $_H{\cal M}$ to the corresponding objects in $\yd $
(see \cite[Proposition 2.4]{bn}).
More precisely, for $B\in {}_H{\cal M}$, ${\bf F}(B)=B$ as a left
$H$-module, and the left $H$-coaction is given by
\begin{equation}\eqlabel{boe}
\l _B(b)=\sum R^2\ot R^1\cd b.
\end{equation}
We observe that $\un {H}^*$ lies in the image of ${\bf F}$, that
is, the relation between action and coaction is given by
\eqref{boe}, which is exactly what we proved in \prref{3.2}.
\end{remarks}

\end{document}